\documentclass[12pt]{amsart}
\usepackage{amscd}
\def\NZQ{\Bbb}               

\def\QQ{{\NZQ Q}}

\def\frk{\frak}               

\def\Phi{{\frk n}}
\def\Phi{{\frk N}}
\def\opn#1#2{\def#1{\operatorname{#2}}} 
\opn\chara{char} \opn\length{\ell} \opn\pd{pd} \opn\rk{rk} \opn\projdim{proj\,dim} \opn\injdim{inj\,dim} \opn\rank{rank} \opn\depth{depth} \opn\grade{grade}
\opn\height{height} \opn\embdim{emb\,dim} \opn\codim{codim}

\opn\Tr{Tr} \opn\bigrank{big\,rank} \opn\superheight{superheight}\opn\lcm{lcm}
\opn\trdeg{tr\,deg}
\opn\reg{reg} \opn\lreg{lreg} \opn\ini{in} \opn\lpd{lpd}
\opn\size{size} \opn\Pf{Pf} \opn\GL{GL} \opn\SL{SL} \opn\mod{mod}
\opn\ord{ord} \opn\Gin{Gin}
\opn\Hilb{Hilb}\opn\adeg{adeg}\opn\std{std}\opn\ip{infpt}
\opn\div{div} \opn\Div{Div} \opn\cl{cl} \opn\Cl{Cl}
\opn\Spec{Spec} \opn\Supp{Supp} \opn\supp{supp} \opn\Sing{Sing} \opn\Ass{Ass} \opn\Min{Min}
\opn\Ann{Ann} \opn\Rad{Rad} \opn\Soc{Soc}
\opn\Syz{Syz} \opn\Im{Im} \opn\Ker{Ker} \opn\Coker{Coker} \opn\Am{Am} \opn\Hom{Hom} \opn\Tor{Tor} \opn\Ext{Ext} \opn\End{End} \opn\Aut{Aut} \opn\id{id}

\opn\nat{nat}
\opn\pff{pf}
\opn\Pf{Pf} \opn\GL{GL} \opn\SL{SL} \opn\mod{mod} \opn\ord{ord} \opn\Gin{Gin} \opn\Hilb{Hilb}
\opn\aff{aff} \opn\con{conv} \opn\relint{relint} \opn\st{st} \opn\lk{lk} \opn\cn{cn} \opn\core{core} \opn\vol{vol} \opn\link{link} \opn\star{star}
\opn\gr{gr}

\def\pot#1#2{#1[\kern-0.28ex[#2]\kern-0.28ex]}

\opn\dirlim{\underrightarrow{\lim}} \opn\inivlim{\underleftarrow{\lim}}
\let\Sect=\bigcap

\let\to=\rightarrow
\let\To=\longrightarrow
\def\Implies{\ifmmode\Longrightarrow \else
        \unskip${}\Longrightarrow{}$\ignorespaces\fi}
\def\implies{\ifmmode\Rightarrow \else
        \unskip${}\Rightarrow{}$\ignorespaces\fi}
\def\iff{\ifmmode\Longleftrightarrow \else
        \unskip${}\Longleftrightarrow{}$\ignorespaces\fi}

\let\:=\colon
\newtheorem{Theorem}{Theorem}[section]
\newtheorem{Lemma}[Theorem]{Lemma}
\newtheorem{Corollary}[Theorem]{Corollary}
\newtheorem{Proposition}[Theorem]{Proposition}

\newtheorem{Example}[Theorem]{Example}

\let\epsilon\varepsilon
\let\phi=\varphi
\let\kappa=\varkappa
\def\qed{\ifhmode\textqed\fi
      \ifmmode\ifinner\quad\qedsymbol\else\dispqed\fi\fi}
\def\textqed{\unskip\nobreak\penalty50
       \hskip2em\hbox{}\nobreak\hfil\qedsymbol
       \parfillskip=0pt \finalhyphendemerits=0}
\def\dispqed{\rlap{\qquad\qedsymbol}}

\opn\dis{dis}
\def\pnt{{\raise0.5mm\hbox{\large\bf.}}}

\opn\Lex{Lex}

\begin{document}

\title{The Strong Lefschetz Property and simple extensions}

\author{J\"urgen Herzog,  Dorin Popescu}


\thanks{The second author was mainly supported by Marie Curie Intra-European
Fellowships MEIF-CT-2003-501046 and partially supported by CNCSIS
and the Ceres program 4-131/2004 of the Romanian Ministery of
Education and Research}

\address{J\"urgen Herzog, Fachbereich Mathematik und
Informatik, Universit\"at Duisburg-Essen, Campus Essen, 45117 Essen,
Germany} \email{juergen.herzog@uni-essen.de}

\address{Dorin Popescu, Institute of Mathematics "Simion Stoilow", University of Bucharest,
P.O.Box 1-764, Bucharest 014700, Romania}
\email{dorin.popescu@imar.ro} \maketitle

\begin{abstract}
Stanley \cite{St} showed that monomial complete intersections have
the strong Lefschetz property. Extending this result  we show that
a simple extension of an Artinian Gorenstein algebra with the
strong Lefschetz property has again the strong Lefschetz property.
\end{abstract}

\section*{Introduction}
Let $K$ be a field,  $A$ be a standard graded Artinian $K$-algebra
and $a\in A$ a homogeneous form   of degree $k$. The element $a$
is called a {\em Lefschetz element} if for all integers $i$ the
$K$-linear map $a\: A_i\to A_{i+k}$ (induced by multiplication
with $a$) has maximal rank. One says that $A$ has  the {\em weak
Lefschetz property} if there exists a Lefschetz element $a\in A$
of degree 1. An element $a\in A_1$ for which all powers $a^r$ are
Lefschetz is called a {\em strong Lefschetz element},  and $A$ is
said to have the {\em strong Lefschetz property } if $A$ admits  a
strong Lefschetz element. Note that the set of Lefschetz elements
$a\in A_1$  form a Zariski open subset of $A_1$. The same holds
true for the set of strong Lefschetz elements.

Assuming that the characteristic of  $K$ is zero and the defining
ideal of $A$ is generated by generic  forms, it is conjectured
that $A$ has the strong Lefschetz property.
 Thus in
particular, $A=K[x_1,\ldots,x_n]/(f_1,\ldots, f_n)$ should have
the strong Lefschetz property for generic forms $f_1,\ldots,f_n$.
Note that such an algebra is an Artinian complete intersection. It
is expected that {\em any} standard graded Artinian  complete
intersection over a base field of characteristic 0  has the strong
Lefschetz property.  Stanley \cite{St} and later J.~Watanabe
\cite{Wa} proved this in case $A$ is a monomial complete
intersection. Stanley used the Hard Lefschetz Theorem to prove
this result, while Watanabe used the representation theory of the
Lie algebra $sl(2)$.

As a main result  of this  paper we prove the following

\medskip
\noindent
 {\bf Theorem:}
\label{main}
 {\em  Let $K$ be a field of characteristic $0$, $A$ be a standard graded Artinian
Gorenstein $K$-algebra having the strong Lefschetz property, and
let $f\in A[x]$ be a monic homogeneous polynomial. Then the
algebra $B=A[x]/(f)$ has the strong Lefschetz property.}

\medskip
The proof only uses techniques from linear algebra. The result
implies in particular Stanley's theorem. More generally it implies
that a complete intersection $K[x_1,\ldots, x_n]/(f_1,\ldots,
f_n)$ with $f_i\in K[x_1,\ldots, x_i]$ for $i=1,\ldots, n$ has the
strong Lefschetz property.

\section{The proof of the main theorem}

Let $A$  be a standard graded $K$-algebra and $I\subset A$ a
graded ideal. For convenience we will say that $a\in A$ is
Lefschetz for $A/I$ if  the residue class $a+I$ is a Lefschetz
element of $A/I$.

In the proof of the main theorem  we shall use the following two
lemmata.

\begin{Lemma}
\label{key1} Let $A$ be a standard graded $K$-algebra, $f,g \in A$
homogeneous elements  which are nonzero divisors on $A$. Then $f$
is Lefschetz  for $A/(g)$  if and only if $g$ is  Lefschetz for
$A/(f)$.
\end{Lemma}

\begin{proof}
Consider the  long exact sequence for Koszul homology (see
\cite[Corollary 1.6.13]{BH})
\[
\begin{CD}
\cdots \to H_1(g;A)\to H_1(f,g;A)\to H_0(g;A)@> f >> H_0(g;A) \to
H_0(f,g;A)\to 0.
\end{CD}
\]
Since $g$ is a non-zerodivisor on $A$ this yields the exact
sequence
\[
\begin{CD}
0\To H_1(f,g;A)\To A/(g)@> f >> A/(g) \to H_0(f,g;A)\to 0.
\end{CD}
\]
Similarly we obtain an exact sequence
\[
\begin{CD}
0\To H_1(f,g;A)\To A/(f)@> g >> A/(f) \to H_0(f,g;A)\to 0.
\end{CD}
\]
Comparing this two exact sequences, the assertion follows.
\end{proof}

\begin{Lemma}
\label{2} Let $K$ be field of characteristic 0, $A$  a
 standard graded Artinian $K$-algebra with strong Lefschetz
property and $f\in A[y]$ a monic homogeneous polynomial. Then for
any strong Lefschetz  element $a\in A_1$ there exists a non-zero
element $c\in K$ such that   $f(a/c)$ is a Lefschetz element of
$A$.
\end{Lemma}
\begin{proof}  Let $f=y^d+a_1 y^{d-1}+\cdots +a_d$, and $s=\max\{i\: A_i\neq 0\}$.
We may assume that $d\leq s$ because otherwise the statement is trivial.
For $c\in K$ we set $f_c=y^d+\sum_{i=1}^d c^i a_i y^{d-i}$. Let $a
\in A_1$  be a strong Lefschetz element. Then  $a^d$ is a
Lefschetz element, that is, for all $i$ the multiplication map
$f_0(a)\: A_i\to A_{i+d}$ has maximal rank.

Fix $i\leq s-d$  and $K$-bases of the nonzero $K$-vector spaces
$A_i$ and $A_{i+d}$, and let $D_c$ be the matrix describing the
$K$-linear map $f_c(a)\: A_i\to A_{i+d}$. Note that the entries of
$D_c$ are polynomial expressions  in $c$ with coefficients in $K$.
Now $P_c(a)$ has maximal rank if and only if one maximal  minor
$M_j(D_c)$  of $D_c$ does not vanish. In particular,
$M_{j_0}(D_0)\neq 0$ for some $j_0$. Since $M_{j_0}(D_c)$ is a
(non-zero) polynomial expression in $c$ with coefficients in $K$,
there exist only finitely many $c\in K$ such that
$M_{j_0}(D_c)=0$. Thus, since $K$ is infinite, we have
$M_{j_0}(D_c)\neq 0$ for infinitely many $c\in K$, and so
$f_c(a)\: A_i\to A_{i+d}$ has maximal rank for infinitely many
$c\in K$. Since $A$ has only finitely many non-zero components, we
can therefore find $c\in K$, $c\neq 0$ such that $f_c(a)$ has
maximal rank for all $i$. Then $a/c\in A_1$ has the desired
property, since $f(a/c)=f_c(a)/c^d$.
\end{proof}

Now we are ready to begin with the proof of the main theorem. Let
$A$ be  a standard graded Artinian Gorenstein $K$-algebra having
the strong Lefschetz property.

In a first step we will prove: suppose $C=A[x]/(x^r)$ has the
strong Lefschetz property for all $r>1$, then $B=A[x]/(f)$ has the
strong Lefschetz property for any monic homogeneous polynomial
$f\in A[x]$.

Let $U_r\subset B_1$ be the Zariski open set of elements $b\in
B_1$ for which $b^r$ is a Lefschetz element. If $U_r\neq
\emptyset$ for all $r\geq 1$, then the finite intersection
$U=\Sect_rU_r$ is non-empty, as well, and any $b\in U$ is then a
strong Lefschetz element. Thus it suffices to show  that for each
$r\geq 1$ there exists an element $b_r\in B_1$ such that $b_r^r$
is a Lefschetz element.

By Lemma \ref{2} we may choose an element $a\in A_1$ such that
$f(a)$ is a Lefschetz element of $A$. It follows that $f(x)$ is
Lefschetz  for $A[x]/(a-x)$.  Thus by Lemma \ref{key1}, the
element $b_1=a-x$ is Lefschetz for $B$.

In case $r>1$, we may view $f(y)$ as a polynomial in $C[y]$ where
$C=A[x]/(x^r)$. By our assumption $C$ has a strong Lefschetz
element. Now Lemma \ref{2} implies that  we can find a strong
Lefschetz element $c\in C_1$ such that $f(c)$ is a Lefschetz
element of $C$. Since the strong Lefschetz elements form a
nonempty Zariski open set in $C_1$, we may assume that
$c=a+\lambda x$ with $a\in A_1$ and $\lambda \in K$, $\lambda\neq
0$. Applying the substitution $x\mapsto b_r=\lambda^{-1}(x-a) $ it
follows that $f(x)$ is Lefschetz
 for $A[x]/b_r^r$. Thus by Lemma \ref{2}, the element
 $b_r^r$ is Lefschetz for $A[x]/(f)$.

\medskip In order to complete the proof of the theorem it remains
to be shown that if $A$ is  a standard graded Artinian Gorenstein
$K$-algebra having the strong Lefschetz property, then
$A[x]/(x^q)$ has the strong Lefschetz property. We use  Lemma
\ref{key1} and show instead that if $a\in A_1$ is a strong
Lefschetz element, then for all $k$ the element $x^q$ is Lefschetz
for $B=A[x]/(a+x)^k$.

\medskip
In $B$ we have
\[
x^k=-\sum_{j=0}^{k-1}\binom{k}{j}a^{k-j}x^j.
\]
By induction on $r$ it follows that
\[
x^r=(-1)^{r-k-1}\sum_{j=0}^{k-1}\binom{r-j-1}{r-k}\binom{r}{j}a^{r-j}x^j\quad
\text{for}  \quad r\geq k.
\]
Note that
\[
\binom{r-j-1}{r-k}\binom{r}{j}=\frac{k}{r-j}\binom{r}{k}\binom{k-1}{j},
\]
so that
\[
x^r=(-1)^{r-k-1}\sum_{j=0}^{k-1}\binom{r}{k}\binom{k-1}{j}\frac{k}{r-j}a^{r-j}x^j
\]
for all $r\geq k$. Thus for all $r\geq 0$ we have
\[
x^r=\sum_{j=0}^{k-1}c_{rj}a^{r-j}x^j
\]
with
\[
c_{rj}= \left\{ \begin{array}{lll} \delta_{rj}, & \mbox{if} & r\leq k-1\\
(-1)^{r-k-1}\binom{r}{k}\binom{k-1}{j}\frac{k}{r-j}, & \mbox{if} &
r\geq k,
\end{array} \right.
\]
where $\delta_{rj}$ denotes the Kronecker symbol.

Now we show that the map $\beta^q_t: B_t\to B_{t+q}$ given by
multiplication with $x^q$ has maximal rank.

We denote by $\alpha^j_i\:A_i\to A_{i+j}$ the $K$-linear map given
by multiplication with $a^j$.   For each element $ux^i\in
A_{t-i}x^i$ we have
\[
x^q(ux^i)=\sum_{j=0}^{k-1}c_{q+i,j}a^{q+i-j}ux^j=\sum_{j=0}^{k-1}c_{q+i,j}\alpha^{q+i-j}_{t-i}(u)x^{j}.
\]
Since for each $j$ the $K$-vectorspace  $B_j$ has the direct sum
decomposition $$B_j=\bigoplus_{i=0}^{k-1} A_{t-i}x^i,$$ the linear
map $\beta^q_t$  can be described by the following block matrix
\[
M=
\begin{pmatrix}
c_{q,0}\alpha^q_t&c_{q+1,0}\alpha^{q+1}_{t-1}&\cdots&c_{q+k-1,0}\alpha^{q+k-1}_{t-k+1}\\
c_{q,1}\alpha^{q-1}_t&c_{q+1,1}\alpha^{q}_{t-1}&\cdots& \cdots\\
\vdots&\vdots&\ddots&\vdots\\
c_{q,k-1}\alpha^{q-k+1}_t&c_{q+1,k-1}\alpha^{q-k+2}_{t-1}&\cdots &
c_{q+k-1, k-1}\alpha^q_{t-k+1}
\end{pmatrix}.
\]
Our aim is to show that $M$ has maximal rank. Assume first that
$q<k$, then
\[
M=\begin{pmatrix} 0 & N\\
\id & *
\end{pmatrix},
\]
where $N=(c_{q+j,i}\alpha^{q+j-i}_{t-j})_{i=0,\ldots, q-1\ \ \atop
j=k-q,\ldots,k-1}$. It follows that $M$ has maximal rank if and
only if $N$ has maximal rank. Thus the general case is treated if
we can prove that for all $q$ the matrix
\[
N=(c_{q+j,i}\alpha^{q+j-i}_{t-j})_{i=0,\ldots, s-1\ \  \atop
j=r-q,\ldots,k-1}\quad \text{with}\quad s=\min\{q,k\}\quad
\text{and}\quad r=\max\{q,k\}
\]
has maximal rank.

We show this by applying certain block row and block column
operations in order to simplify the matrix without changing its
rank. The kind of operations we will apply are the following:
\begin{enumerate}
\item[(i)] multiplication of a block row or a block column of $N$
with a non-zero rational number;

\item [(ii)]  for $d\in \QQ$ , $d\neq 0$ and $j<i$ compose
$d\alpha^{i-j}_{t-i}$ with each block $c_{jl}\alpha^{j-l}_{t-j}$
of the $j$th block column of $N$ to obtain a $j$th block column
whose blocks are $dc_{jl}\alpha^{j-l}_{t-j}\circ
\alpha^{i-j}_{t-i}=dc_{jl}\alpha^{i-l}_{t-i}$, and subtract this
new block column from the $i$th block column to obtain the new
$i$th block column whose blocks are
$$(c_{il}-dc_{jl})\alpha^{i-l}_{t-i},\quad l=0,\cdots, s-1.$$
\end{enumerate}

 These operations only change the coefficients
$c_{ji}$ of the block entries, that is, the matrix
$N=(c_{q+j,i}\alpha^{q+j-i}_{t-j})_{i=0,\ldots, s-1\ \  \atop
j=r-q,\ldots,k-1}$ will be transformed into a matrix of the form
$N'=(c'_{q+j,i}\alpha^{q+j-i}_{t-j})_{i=0,\ldots, s-1\ \  \atop
j=r-q,\ldots,k-1}$ with certain new coefficients $c'_{q+j,i}\in
\QQ$.

\medskip
Consider the  ``coefficient matrix" $L=(c_{q+j,i})_{i=0,\ldots,
s-1\ \ \atop j=r-q,\ldots,k-1}$ of $N$. Then the coefficient
matrix $L'$ of $N'$ is obtained from $L$ by the following row and
column operations:
\begin{enumerate}
\item[(i)] multiplication or division of  a row or a column with a
non-zero rational number;

\item[(ii)] subtraction of a multiple of the $j$th column from the
$i$th column where $j<i$.
\end{enumerate}

Next we intend  to show that by these operations $L$ can be
transformed into a matrix $L'$ such that all entries of $L'$ on
the anti-diagonal are non-zero, while the entries below the
anti-diagonal are all zero.

We first simplify $L$ by dividing  each $j$th column by
$(-1)^{j-k-1}{j \choose k}k$ and each $i$th row by ${k-1\choose
i}$. The result of these operations is the matrix
\[
\begin{pmatrix}
1/r & 1/(r+1) &  \cdots & 1/(r+s-1)\\
1/(r-1) & 1/r &\cdots & 1/(r+s-2)\\
\vdots & \vdots&  \ddots & \vdots\\
1/(r-s+1)& 1/(r-s+2) & \cdots & 1/r
\end{pmatrix},
\]
which we again denote by $L$.

We will use the following simple fact from linear algebra: suppose
$F=(f_{ij})_{i,j=1,\ldots,n}$ is an $n\times n$-matrix with
coefficients in a field $K$. Then the following conditions are
equivalent:
\begin{enumerate}
\item[(a)] the matrix $F$ can be transformed by  operations of
type (ii) into a matrix $F'$ with $f'_{ij}\neq 0$ for $i+j=n+1$,
and $f'_{ij}=0$ for  $i+j>n+1$;

\item[(b)]  $\det(F_i)\neq 0$ for $i=1,\ldots,n$ where
$F_i=(f_{kl})_{k=i,\ldots,n\atop l=1,\ldots i}$.
\end{enumerate}

Indeed, it is clear that (a)\implies (b). Conversely, assuming (b)
we have $\det(F_n)=f_{n1}\neq 0$. Thus by subtracting suitable
multiples of  the first column from the other columns we obtain a
matrix $G=(g_{ij})$ with $g_{ni}= 0$ for $i=2,\ldots, n$, and such
that $\det(F_i)=\det(G_i)$ for $i=1,\ldots, n$, where
$G_i=(g_{kl})_{k=i,\ldots,n\atop l=1,\ldots i}$. Applying the
induction hypothesis to the matrix
$G'=(g_{kl})_{k=1,\ldots,n-1\atop l=2,\ldots, n}$, the assertion
follows.

\medskip
Applying this result from linear algebra, we see that $L$ can be
transformed by operations of type (ii) into the matrix $L'$ of the
desired form if for all integers $0\leq t<s$ the  matrices of the
shape
\[
S=(1/(r-i+j))_{i,j=0,\ldots,t}
\]
are non-singular. It is  an easy exercise in linear algebra to
show that this is indeed the case.

After all these operations our matrix $N$ is transformed into the
matrix $N'$ whose anti-diagonal has the block entries
\[
c_1' \alpha^{r-s+1}_{t+q-r}, c_2'\alpha^{r-s+3}_{t+q-r-1},\ldots,
c_{s-1}'\alpha^{q+k-1}_{t-k+1}
\]
with non-zero rational coefficients $c_i'$, and whose block
entries below the anti-diagonal are all zero.

We will show that for $i=0,\ldots,s-1$ either all
$\alpha^{r-s+2i+1}_{t+q-r-i}$ are injective maps, or else all
$\alpha^{r-s+2i+1}_{t+q-r-i}$ are surjective maps. Then clearly
$N'$ has maximal rank, and consequently $N$ has maximal rank.

For all integers $i$ and $j$ with $0\leq i<j$ the maps
\[
\alpha_i^{j-i}\: A_i\To A_j
\]
have maximal rank, by assumption. In particular, $\alpha_i^{j-i}$ is
injective if $\dim A_i\leq \dim A_j$ and  surjective if   $\dim A_i\geq \dim A_j$.

Let $\sigma=\max\{i\:A_i\neq 0\}$. Then, since $A$ is Gorenstein,
the Hilbert function of $A$ is symmetric (see e.g. \cite[Corollary
4.4.6, Remark 4.4.7]{BH}), that is,
\[
\dim A_i=\dim A_{\sigma-i} \quad \text{for all}\quad i,
\]
and since  $A$ has the weak Lefschetz property ($A$ even has the
strong Lefschetz property), the Hilbert function of $A$ is
unimodal (see e.g. \cite[Remark 3.3]{HMNW}).  It then follows that
\[
\dim A_i\leq   \dim A_j \quad\text{if and only if}\quad i\leq \sigma-j.
\]
Thus we conclude that
\[
 \alpha_i^{j-i} \quad\text{is}\quad  \left\{ \begin{array}{lll} \text{injective} & \text{if} & i\leq \sigma-j,\\
 \text{surjective} & \mbox{if} & i\geq \sigma-j.
 \end{array} \right.
 \]
 Thus in case of the maps $\alpha^{r-s+2i+1}_{t+q-r-i}$,
 we have to compare the size of the numbers
 $t+q-r-i$ and  $\sigma -[(r-s+2i+1)+(t+q-r-i)]=\sigma-t-q+s-i-1$. Since it does not depend on $i$  which of the
 two numbers is less than or equal to other, it follows  $\alpha^{r-s+2i+1}_{t+q-r-i}$ is injective for all $i$,
 or   $\alpha^{r-s+2i+1}_{t+q-r-i}$  is surjective for all $i$, as
 desired.
\section{Some comments}

As an immediate consequence of our main theorem   we obtain

\begin{Corollary}
\label{cor1} Let $K$ be field of characteristic 0, and  $A$ be an
Artinian Gorenstein $K$-algebra. For $i=1,\ldots, n$ let $f_i\in
A[x_1,\ldots,x_i]$ be a homogeneous polynomial which is monic in
$x_i$. Then  the $K$-algebra
\[
A[x_1,\ldots, x_n]/(f_1,\ldots, f_n)
\]
has the strong Lefschetz property.
\end{Corollary}

 The result implies in particular that $K[x_1,\ldots,
x_n]/(f_1,\ldots, f_n)$ has the strong Lefschetz property,  if for
$i=1,\ldots, n$, $f_i\in K[x_1,\ldots,x_i]$ is a homogeneous and
monic polynomial in $x_i$. In the special case that
$f_i=x_i^{a_i}$ for $i=1,\ldots,n$, we obtain the theorem of
Stanley \cite{St}. The slightly more general result with the $f_i$
as described before, can also be deduced directly from Stanley's
theorem using the following result of Wiebe \cite[Proposition
2.9]{Wi}: let $I\subset K[x_1,\ldots, x_n]$ be a graded ideal, and
assume that $K[x_1,\ldots, x_n]/\ini(I)$ has the strong Lefschetz
property, where $\ini(I)$ is the initial ideal with respect to
some term order. Then $K[x_1,\ldots, x_n]/I$  has the strong
Lefschetz property.

In the above situation we have $\ini(f_i)=x_i^{\deg f_i}$ for $i=1,\ldots,n$, if we choose the
lexicographical order induced by $x_n>x_{n-1}>\cdots > x_1$. Since the initial terms of the generators
form   a regular sequence
it follows that $$\ini(I)=(\ini(f_1),\ldots, \ini(f_n))=(x_1^{a_1},\ldots, x_n^{a_n}).$$

\bigskip
In case the  $K$-algebra $A$ is not Gorenstein,  the proof of our
main theorem yields the following weaker result.

\begin{Proposition} \label{weak}
 Let $K$ be a field, $A$  a standard graded Artinian
$K$-algebra having the strong Lefschetz property, and let $f\in
A[x]$ be a monic homogeneous polynomial. Then the algebra
$B=A[x]/(f)$ has the weak Lefschetz property.
\end{Proposition}

\begin{proof} Recall the following step in the proof of the main theorem: by Lemma \ref{2} we may choose an element $a\in A_1$ such that
$f(a)$ is a Lefschetz element of $A$. It follows that $f(x)$ is
Lefschetz  for $A[x]/(a-x)$.  Thus by Lemma \ref{key1}, $b=a-x$ is
Lefschetz for $B$.
\end{proof}

Analyzing the arguments in the proof of our main theorem we see
that all  results remain valid if the  characteristic of the base
field is large enough. More precisely we have

\begin{Corollary}
\label{cor2} Let $K$ be a field and $A$  an Artinian Gorenstein
$K$-algebra with socle degree $\sigma=\max\{t\: A_t\not=0\}$ and
multiplicity  $e(A)=\sum_{t=0}^\sigma\dim_KA_t$. Let $f\in A[x]$
be a homogeneous monic polynomial of degree $q$.  Then
$B=A[x]/(f)$ has the strong Lefschetz property if
\[
\chara K\geq \left\{ \begin{array}{lll} 2q+\sigma -1,& \text{and $f=x^q$,} & \\
\max\{e(A), 2q+\sigma -1\}, & \text{otherwise}.
 \end{array} \right.
 \]
\end{Corollary}

\begin{proof} In case $f=x^q$ we must make sure that all the
binomials in the expression
$x^r=(-1)^{r-k-1}\sum_{j=0}^{k-1}\binom{r-j-1}{r-k}\binom{r}{j}a^{r-j}x^j$
are units in the field $K$, and this must be satisfied for all
$r=q+k$ where  are less than or equal the socle degree of
$A[x]/(x^q)$. Since the socle degree of $A[x]/(x^q))$ is equal to
$q+\sigma-1$, we therefore need that $\chara K$ does not divide
any prime number $\leq 2q +\sigma -1$.

In the general case we had to apply Lemma \ref{2}. For the proof
of this lemma it was necessary that the field  $K$ has enough
elements, so that for all the polynomials in $c$ defined by the
maximal minors considered in the proof we find a common element
$c\in K$ for which these polynomials do not vanish. This is
possible if $\chara K>e(A)$.
\end{proof}

We conclude this note with the following

\begin{Example}
\label{gegen} {\em If $A$ has the  strong Lefschetz property and
$f\in A$ is a generic form. One expect that $A/(f)$ has again the
strong Lefschetz property.  However, in general this is no the
case. Indeed, let
$A=K[x_1,\ldots,x_5]/(x_1^4,x_2^4,x_3^4,x_4^4,x_5^2)$. Then $A$
has the strong Lefschetz property. Let $f\in A$ be a generic form
of degree 8 and set $B=A/(f)$. (We use the ``Randomized" command
of CoCoA to produce generic forms.)

The Hilbert series of $B$ is given by
$$\text{Hilb}_B(t)=1+5t+14t^2+30t^3+51t^4+71t^5+84t^6+84t^7+70t^8+46t^9+16t^{10}.$$
Let $b\in B$ be a generic linear form, and set $C=B/(b^9)$. Then
\[
\text{Hilb}_{C}(t)=1+5t+14t^2+30t^3+51t^4+71t^5+84t
^6+84t^7+70t^8+45t^9+12t^{10}.
\]
 It follows that the  map $ B_1\xrightarrow{b^9} B_{10}$ is not
surjective but also not injective because
$\dim_KB_1+\dim_KC_{10}=5+12>16=\dim_KB_{10}$. Thus $B$ does not
have the strong Lefschetz property.

On the other hand it can be checked that $B$   has the maximal
rank property, that is, any generic form in $B$ has maximal rank.
Such an example seems to be new, see \cite{MR}. }
\end{Example}


\begin{thebibliography}{1}


\bibitem{BH}  W.\ Bruns,\ J.\ Herzog, {\sl Cohen-Macaulay rings}, Revised
Edition,  Cambridge, 1996.

\bibitem{HMNW} T.\ Harima,\ J.\ Migliore,\ U.\ Nagel,\ J.\
Watanabe, The weak and strong Lefschetz properties for artinian
$K$-algebras, arXiv:math.AC/0208201.


\bibitem{MR} J.\ Migliore,\ R.M.\ Miro-Roig, Ideals of generic forms and the ubiquity
of the weak Lefschetz property, arXiv:math.AC/0205133.

\bibitem{St} R.\ Stanley, Weyl groups, the hard Lefschetz theorem,
and the Sperner property, SIAM J. Algebraic Discrete Methods 1
(1980), 168-184.


\bibitem{Wa} J.\ Watanabe, The Dilworth number of Artinian rings and
finite posets with rank function, Commutative Algebra and
Combinatorics, Advanced Studies in Pure Math.,Vol {\bf 11},
Kinokuniya Co. North Holland, Amsterdam, (1987), 303-312.

\bibitem{Wi} A.\ Wiebe, The Lefschetz property for componentwise
linear ideals and Gotzmann ideals, arXiv:math.AC/0307223.
\end{thebibliography}
\end{document}